\documentclass[runningheads]{llncs}
\usepackage{graphicx}

\usepackage{amsmath} 
\usepackage{amsfonts}
\usepackage{amssymb}
\usepackage[colorlinks=true,breaklinks,pdfborder={0 0 0},linkcolor=black,citecolor=black,urlcolor=blue]{hyperref}

\newcommand{\ZFC}{\ensuremath{\mathsf{ZFC}} }
\newcommand{\ZF}{\ensuremath{\mathsf{ZF}} }

\newcommand{\HOD}{\ensuremath{\mathsf{HOD}} }

\newcommand{\AD}{\ensuremath{\mathsf{AD}} }

\newcommand{\BS}{{}^\omega\omega}

\newcommand{\PI}{\boldsymbol\Pi}
\newcommand{\SIGMA}{\boldsymbol\Sigma}

\newcommand{\cM}{\mathcal{M}}

\newcommand{\bR}{\mathbb{R}}

\newcommand{\Pot}{\mathcal{P}}

\begin{document}
\title{Determinacy Axioms and Large Cardinals\thanks{Supported by L'OR\'{E}AL Austria, in collaboration with the Austrian UNESCO Commission and in cooperation with the Austrian Academy of Sciences - Fellowship \emph{Determinacy and Large Cardinals} and the Austrian Science Fund (FWF) under Elise Richter grant number V844, international project number I6087, and START grant number Y1498.}}
%
%
\author{Sandra Müller\inst{1}\orcidID{0000-0002-7224-187X}}
\authorrunning{S. Müller}
%
\institute{Institut f\"ur Diskrete Mathematik und Geometrie, TU Wien,\\ Wiedner Hauptstra{\ss}e 8-10/104, 1040 Wien, Austria\\
\email{sandra.mueller@tuwien.ac.at}\\
\url{https://dmg.tuwien.ac.at/sandramueller/}}
\maketitle              
\begin{abstract}

The study of inner models was initiated by Gödel’s analysis of the constructible universe. Later, the study of canonical inner models with large cardinals, e.g., measurable cardinals, strong cardinals or Woodin cardinals, was pioneered by Jensen, Mitchell, Steel, and others. Around the same time, the study of infinite two-player games was driven forward by Martin’s proof of analytic determinacy from a measurable cardinal, Borel determinacy from ZFC, and Martin and Steel’s proof of levels of projective determinacy from Woodin cardinals with a measurable cardinal on top. First Woodin and later Neeman improved the result in the projective hierarchy by showing that in fact the existence of a countable iterable model, a mouse, with Woodin cardinals and a top measure suffices to prove determinacy in the projective hierarchy. This opened up the possibility for an optimal result stating the equivalence between local determinacy hypotheses and the existence of mice in the projective hierarchy. This article outlines the main concepts and results connecting determinacy hypotheses with the existence of mice with large cardinals as well as recent progress in the area.

\end{abstract}
\section{Introduction} 

The standard axioms of set theory, Zermelo-Fraenkel set theory with Choice ($\ensuremath{\mathsf{ZFC}}$), do not suffice to answer all questions in mathematics. While this follows abstractly from Kurt Gödel's famous incompleteness theorems, we nowadays know numerous concrete examples for such questions.
A large number of problems in set theory, for example, regularity properties such as Lebesgue measurability and the Baire property are not decided -- for even rather simple (for example, projective) sets of reals -- by $\ensuremath{\mathsf{ZFC}}$. Even many problems outside of set theory have been showed to be unsolvable, meaning neither their truth nor their failure can be proven from $\ensuremath{\mathsf{ZFC}}$. This includes the Whitehead Problem (group theory, \cite{Sh74}), the Borel Conjecture (measure theory, \cite{La76}), Kaplansky's Conjecture on Banach algebras (analysis, \cite{DW87}), and the Brown-Douglas-Fillmore Problem (operator algebras, \cite{Fa11}). A major part of set theory is devoted to attacking this problem by studying various extensions of $\ZFC$ and their properties. One of the main goals of current research in set theory is to identify \emph{the ``right'' axioms for mathematics} that settle these problems. This, in part philosophical, problem is attacked with technical mathematical methods by analyzing various extensions of $\ZFC$ and their properties.
{\bf Determinacy assumptions} are canonical extensions of $\ensuremath{\mathsf{ZFC}}$ that postulate the existence of winning strategies in natural two-player games. Such assumptions are known to imply regularity properties, and 
enhance sets of real numbers with a great deal of canonical structure. Other natural and well-studied extensions of $\ZFC$ are given by the hierarchy of {\bf large cardinal axioms}. Determinacy assumptions, large cardinal axioms, and their consequences are widely used and have many fruitful implications in set theory and even in other areas of mathematics such as algebraic topology \cite{CSS05}, topology \cite{Ny80,Fl82,CMM1}, algebra \cite{EM02}, and operator algebras \cite{Fa11}. Many applications, in particular, proofs of consistency strength lower bounds, exploit the interplay of large cardinals and determinacy axioms. Thus, understanding the connections between determinacy assumptions and the hierarchy of large cardinals is vital to \textbf{answer questions left open by $\ZFC$ itself}. The results outlined in this article are closely related to this overall goal.


To explore the connections between large cardinals and determinacy at higher levels, the study of other hierarchies, for example, with more complex inner models called \textbf{hybrid mice}, has been very fruitful.  \textbf{Translation procedures} are needed to translate these hybrid models, whose strength comes from descriptive set theoretic features, back to standard inner models while making use of their hybrid nature to obtain stronger large cardinals in the translated model. They are therefore a key method \textbf{connecting descriptive set theory with inner model theory}.
One of the results surveyed in this article is a new translation procedure extending work of Sargsyan \cite{Sa17}, Steel \cite{St08dermod}, and Zhu \cite{Zh15}. This new translation procedure yields a countably iterable inner model with a cardinal $\lambda$ that is both a limit of Woodin cardinals and a limit of strong cardinals \cite{Mu21}. So it improves Sargsyan's construction in \cite{Sa17} in two ways: It can be used to obtain infinitely many instead of finitely many strong cardinals and the models it yields are countably iterable -- a crucial property of mice. This translation procedure can be applied to prove a \textbf{conjecture of Sargsyan} on the consistency strength of the Axiom of Determinacy when all sets are universally Baire \cite{Mu21}, a central and widely used property of sets of reals introduced implicitly in \cite{SchVa83} and explicitly in \cite{FMW92}. In fact, the new translation procedure can be applied in a much broader context. Moreover, it provides the basis for translation procedures resulting in more complex patterns of strong cardinals, for example, a strong cardinal that is a limit of strong cardinals.

Recent seminal results of Sargsyan and Trang \cite{SaTrC,SaTrB,SaTr}, see also the review \cite{Mu20rev}, as well as Larson and Sargsyan \cite{LaSa21,Sa21} suggest that we are at a \textit{turning point} in the search for natural constructions of canonical models with a Woodin limit of Woodin cardinals and thereby for proving better lower bounds for natural set theoretic hypotheses.

\section{Determinacy for games of length $\omega$ and large cardinals}

In 1953, Gale and Stewart \cite{GS53} developed a basic theory of infinite games. For notational simplicity, we identify reals in $\bR$ with $\omega$-sequences of natural numbers in $\BS$. Gale and Stewart considered, for every set of reals $A$, a two-player game
$G(A)$ of length $\omega$, where player $\mathrm{I}$ and player
$\mathrm{II}$ alternate playing natural numbers $n_0, n_1, \dots$, as follows:
\[ \begin{array}{c|ccccc} \mathrm{I} & n_0 & & n_2 & &\hdots \\ \hline
    \mathrm{II} & & n_1 & & n_3 & \hdots 
   \end{array} \]

They defined that
player $\mathrm{I}$ wins the game $G(A)$ if and only if the sequence $x = (n_0,n_1,\dots)$ of natural numbers produced during a run of the game $G(A)$ is
an element of $A$; otherwise, player $\mathrm{II}$ wins. We call $A$ the payoff set of $G(A)$. The
game $G(A)$ (or the set $A$ itself) is called determined if and only if
one of the two players has a winning strategy, meaning that there is a method by
which they can win in the game described above, no matter what their opponent does. The \textbf{Axiom of Determinacy} ($\AD$) is the statement that all sets of reals are determined.

Already in \cite{GS53}, the authors were able to prove that every open
and every closed set of reals is determined under $\ZFC$. But
they also proved that determinacy for all sets of reals
contradicts the Axiom of Choice. This leads to the natural question as to
how the picture looks for definable sets of reals which are more
complicated than open and closed sets. After some partial results by Wolfe \cite{Wo55}
and Davis \cite{Da64},
Martin was able to prove in 1975 \cite{Ma75} that
every Borel set of reals is determined (again using $\ZFC$).

In the meantime, the development of so-called \textbf{large cardinal
  axioms} was proceeding in set theory, and Solovay was able to prove regularity properties, a known consequence of determinacy, for a specific pointclass, assuming the existence of a measurable cardinal, instead of a determinacy axiom. 
  Finally, Martin was able to prove a direct connection between large cardinals and determinacy axioms: He showed, in 1970, that the
existence of a measurable cardinal implies determinacy for every
analytic set of reals \cite{Ma70}. Eight years later, Harrington  established that this result is, in some sense, optimal, by proving that determinacy for all analytic sets of reals implies that $0^\#$, a countable active iterable canonical inner model which can be obtained from a measurable cardinal, exists \cite{Ha78}. Here, an {\bf iterable canonical inner model}, or {\bf mouse}, is a fine structural model that is, in some sense, iterable. This notion goes back to Jensen \cite{Je72}. Together with Martin's argument mentioned above, this yields an equivalence between the two statements. The construction of such canonical inner models and their connection with determinacy was later extended in work of Dodd, Jensen, Martin, Mitchell, Neeman, Schimmerling, Schindler, Solovay, Steel, Woodin, Zeman, and others (see, e.g., \cite{MaSt89,MS94,SchStZe02,St10,Ze02}; 
see the preface of \cite{MSW20} or Larson's history of determinacy \cite{La} for a more detailed overview). In the projective hierarchy, this led to the following fundamental theorem. Here, $M_0^\#(x)$ denotes $x^\#$, a version of $0^\#$ relativized to a real $x$, and $M_n^\#(x)$ denotes a minimal countable active mouse with $n$ Woodin cardinals constructed above $x$. 

\begin{theorem}[Harrington, Martin, Neeman, Woodin \cite{Ha78,Ma70,Ne95,Ne02,MSW20}]
  Let $n$ be a natural number. Then the following are equivalent:
  \begin{enumerate}
  \item All $\PI^1_{n+1}$ sets are determined, and
  \item for all $x \in \BS$, $M_n^\#(x)$ exists and is $\omega_1$-iterable.
  \end{enumerate}
\end{theorem}

The proof that the determinacy of sets in the projective hierarchy implies the existence of mice with finitely many Woodin cardinals in this exact level-by-level correspondence first appeared in \cite{Uh16,MSW20} and is originally due to Woodin. As shown in \cite{Uh16}, the underlying methods can be used to obtain similar results for certain hybrid mice in the $L(\bR)$-hierarchy.
These tight connections are, at first, very surprising, as they show that two ostensibly completely different notions, from distinct areas of set theory -- determinacy from descriptive set theory, and inner models with large cardinals from inner model theory -- are, in fact, the same.

\section{Determinacy for games longer than $\omega$}

It turns out that the correspondence between determinacy and inner models with large cardinals does not stop at games of length $\omega$. For every ordinal $\alpha$ and set $A \subseteq {}^\alpha\omega$, we can define a game $G(A)$ of length $\alpha$ with payoff set $A$, as follows:
\[ \begin{array}{c|cccccccc} \mathrm{I} & n_0 & & n_2 & &\hdots &
    n_\omega & & \hdots \\ \hline
    \mathrm{II} & & n_1 & & n_3 & \hdots & & n_{\omega+1} & \hdots 
   \end{array} \]

 The players alternate playing natural numbers $n_i$ for $i<\alpha$, and we again say that player $\mathrm{I}$ wins the game if and only if the sequence $x = (n_0,n_1,\dots)$ of length $\alpha$ they produce is an element of $A$; otherwise, player $\mathrm{II}$ wins. In landmark results, Neeman \cite{Ne04} developed powerful techniques to prove the determinacy of projective games longer than $\omega$ from large cardinals. A first step in this direction is, for example, the following result:

\begin{theorem}[Neeman, \cite{Ne04}]
  Let $n \in \omega$ and suppose that $M_{\omega+n}^\#(x)$ exists for all reals $x \in \BS$. Then all games of length $\omega^2$ with $\PI^1_{n+1}$ payoff are determined.
\end{theorem}

This result in fact holds for games of fixed length $\alpha$, for all countable ordinals $\alpha$, instead of games of length $\omega^2$. The following theorem complements Neeman's results for projective games of length $\omega^2$:

\begin{theorem}[Aguilera, Müller, \cite{AgMu19,Mu19}]
  Let $n$ be a natural number and suppose that all games of length $\omega^2$ with $\PI^1_{n+1}$ payoff are determined. Then, for every $x \in \BS$, there is a model $\cM$ of $\ZFC$, with $\omega+n$ Woodin cardinals, such that $x \in \cM$.
\end{theorem}

At this level, the interplay of determinacy and large cardinals is already understood quite well (see also \cite{AgMu,AMS}). For games of length $\omega^\alpha$ with analytic payoff, for countable ordinals $\alpha$, similar results have previously been established by Trang \cite{Tr13}, building on unpublished results of Woodin, using canonical models of determinacy with a generalized Solovay measure. The Solovay measure is also called a \textbf{supercompact measure for $\omega_1$} as it witnesses a degree of supercompactness for $\omega_1$. 

%
When considering much stronger notions of determinacy, the picture is less clear. For example, it was already shown by Mycielski in 1964 that determinacy for all games of length $\omega_1$ is inconsistent with Zermelo–Fraenkel set theory ($\ZF$). Nevertheless, there are subclasses of games of length $\omega_1$ that are still known to be determined under large cardinal assumptions.

An intermediate step are games that do not have a fixed countable length but still end after countably many rounds. In 2004, Neeman showed in groundbreaking work, from large cardinals, that certain games that are not of fixed countable length are still determined. These so-called {\bf games of continuously coded length}, which go back to Steel \cite{St88}, are defined as follows: For any set $A \subset (\BS)^{{<}\omega_1}$ and partial function $\nu \colon \BS \rightharpoonup \omega$, the game $G_{\text{cont}}(\nu,A)$ is given by the following rules:
\[ \begin{array}{c|cccccccccc} \mathrm{I} & y_0(0) & & y_0(2) & &\hdots &
    y_\alpha(0) & & y_\alpha(2) & & \hdots \\ \hline
    \mathrm{II} & & y_0(1) & & y_0(3) & \hdots & & y_\alpha(1) & & y_\alpha(3) & \hdots 
   \end{array} \]

We canonically identify segments of the game of length $\omega$ as mega-rounds, and let $y_\alpha$ denote the real that the two players together produce in mega-round $\alpha$. If $\nu(y_\alpha)$ is undefined, the game ends, and player $\mathrm{I}$ wins if and only if $\langle y_\xi \mid \xi \leq \alpha \rangle \in A$. Otherwise, let $n_\alpha =  \nu(y_\alpha)$. Then the game ends if $n_\alpha \in \{n_\xi \mid \xi < \alpha\}$, and again, player $\mathrm{I}$ wins if and only if $\langle y_\xi \mid \xi \leq \alpha \rangle \in A$. If neither of these alternatives hold, the game continues.

\begin{theorem}[Neeman, \cite{Ne04}]\label{thm:contcoded}
  Suppose there is an iterable proper class model $M$, with a Woodin cardinal $\delta$ and a cardinal $\kappa < \delta$ that is $(\delta+1)$-strong in $M$, such that $V^M_{\delta+1}$ is countable in $V$. Then the game $G_{\operatorname{cont}}(\nu,A)$ is determined for every $\nu$ in the class $\SIGMA^0_2$ and every $A$ that is ${<}\omega^2-\PI^1_1$ in the codes.
\end{theorem}

Here, being $\Gamma$ in the codes for a pointclass $\Gamma$ and a set $A \subset (\BS)^{{<}\omega_1}$ is defined via a natural coding of elements of $A$ as reals; $A$ is \emph{$\Gamma$ in the codes} if the set of codes of elements of $A$ belongs to $\Gamma$.
It is not known whether Theorem \ref{thm:contcoded} is optimal, but results of Neeman and Steel \cite{Ne99} show that it cannot be very far away from optimal.
I conjecture that it is indeed optimal in the following sense:

\begin{conjecture}\label{conj:contcoded}
  Suppose the game $G_{\operatorname{cont}}(\nu,A)$ is determined for every $\nu$ in the class $\SIGMA^0_2$ and every $A$ that is ${<}\omega^2-\PI^1_1$ in the codes. Then there is a model of $\ZFC$ with a Woodin cardinal $\delta$ and a cardinal $\kappa < \delta$ that is $(\delta+1)$-strong.
\end{conjecture}

A similar conjecture at a higher level is moving toward a Holy Grail of current inner model theory. More precisely, it concers the aim to prove the existence of an inner model with a Woodin cardinal that is a limit of Woodin cardinals from the determinacy of certain long games. The natural games to consider at this level have length $\omega_1$ and their payoff set is ordinal definable using reals as parameters. The converse was shown by Woodin, using results of Neeman \cite{Ne04} and ideas going back to Kechris and Solovay \cite{KS85}.

\begin{theorem}[Neeman, Woodin, \cite{Ne04}]\label{thm:WoodinODGamesOmega1}
  Suppose there is an iterable proper class model with a Woodin cardinal that is a limit of Woodin cardinals and countable in $V$. Then there is a model of $\ZFC$ in which all ordinal definable games of length $\omega_1$ on natural numbers with real parameters are determined.
\end{theorem}

In fact, Woodin showed that determinacy of these games of length $\omega_1$ is equiconsistent with a seemingly weaker statement: determinacy of certain games that are {\bf constructibly uncountable in the play}. These games are defined as follows: For a payoff set $A \subset (\BS)^{{<}\omega_1}$, players $\mathrm{I}$ and $\mathrm{II}$ alternate playing natural numbers to produce reals $y_\alpha$.
\[ \begin{array}{c|cccccccccc} \mathrm{I} & y_0(0) & & y_0(2) & &\hdots &
    y_\alpha(0) & & y_\alpha(2) & & \hdots \\ \hline
    \mathrm{II} & & y_0(1) & & y_0(3) & \hdots & & y_\alpha(1) & & y_\alpha(3) & \hdots 
   \end{array} \]

The game ends when its length reaches the first ordinal $\gamma$ which is uncountable in $L[y_\alpha \mid \alpha < \gamma]$, and player $\mathrm{I}$ wins if and only if $\langle y_\alpha \mid \alpha < \gamma \rangle \in A$. Here $L[y_\alpha \mid \alpha < \gamma]$ denotes Gödel's Constructible Universe $L$ relative to $(y_\alpha \mid \alpha < \gamma)$. In this case $\gamma = \omega_1^{L[y_\alpha \mid \alpha < \gamma]}$, so it makes sense to say that the game {\bf ends at $\omega_1$ in $L$ of the play}. We technically define that $\mathrm{II}$ wins if the game lasts $\omega_1$ (in $V$) rounds, but mild large cardinal assumptions yield an ordinal $\gamma$, as above, that is countable in $V$. Neeman proved that, for sufficiently definable payoff sets $A$, these games are determined, via a sophisticated extension of the methods used in the proof of Theorem \ref{thm:contcoded}.

\begin{theorem}[Neeman, \cite{Ne04}]\label{thm:NeemanGamesUnctbleInL}
  Suppose there is an iterable proper class model with a Woodin cardinal that is a limit of Woodin cardinals and countable in $V$. Then all games ending at $\omega_1$ in $L$ of the play with payoff sets that are $\Game({<}\omega^2-\Pi^1_1)$ in the codes are determined.
\end{theorem}

Here, $\Game$ denotes the game quantifier for games of length $\omega$. In \cite{Ne02wlw}, Neeman showed the consistency of the hypotheses of Theorems \ref{thm:WoodinODGamesOmega1} and \ref{thm:NeemanGamesUnctbleInL} from large cardinals. In light of Theorem \ref{thm:NeemanGamesUnctbleInL}, Theorem \ref{thm:WoodinODGamesOmega1} is a consequence of the following result of Woodin's:

\begin{theorem}[Woodin, \cite{Ne04}]\label{thm:WoodinDetEquiconsistency}
The following theories are equiconsistent:
\begin{enumerate}
\item $\ZFC \, +$ all ordinal definable games of length $\omega_1$ on natural numbers with real parameters are determined.
\item $\ZFC \, +$ all games ending at $\omega_1$ in $L$ of the play with payoff sets that are $\Game({<}\omega^2-\Pi^1_1)$ in the codes are determined.
\end{enumerate}
\end{theorem}

I conjecture that Theorem \ref{thm:WoodinODGamesOmega1} is optimal, in the following sense:

\begin{conjecture}\label{conj:wlw}
  Suppose all ordinal definable games of length $\omega_1$ on natural numbers with real parameters are determined. Then there is a model of $\ZFC$ with a Woodin cardinal that is a limit of Woodin cardinals.
\end{conjecture}

This would be the first correspondence between a natural determinacy notion and large cardinals at the level of a Woodin cardinal that is a limit of Woodin cardinals. It cannot be achieved using current methods such as the core model induction technique due to Woodin (see, for example, the review \cite{Mu20rev}), which Sargsyan and Trang \cite{SaTrC,SaTr,SaTrB} have recently shown runs into serious issues before reaching this level. In addition, by recent results of Larson and Sargsyan \cite{LaSa21,Sa21}, also the well-known liberal $K^c$ construction in \cite{ANS01,JSSS09} can fail if there is a Woodin cardinal that is a limit of Woodin cardinals.

 


Therefore, understanding the large cardinal strength of the determinacy of such uncountable games might shed light on how to canonically obtain inner models with a Woodin cardinal that is a limit of Woodin cardinals.

\section{Strong models of determinacy for games of length $\omega$}

Another approach to strengthen determinacy is to keep playing games of length $\omega$ and impose additional structural properties on the model. Examples of such structural properties are ``$\theta_0 < \Theta$,'' ``$\Theta$ is regular,'' or the Largest Suslin Axiom, see, for example, \cite{St08dermod,Sa15,SaTr}. Here $\Theta$ is given by \[ \Theta = \sup\{ \alpha \mid \text{there is a surjection } f \colon \bR \rightarrow \alpha \} \] and we write
$\theta_0$ for the least ordinal $\alpha$ such that there is no surjection of $\bR$ onto $\alpha$ which is ordinal definable from a real. While in models of the Axiom of Choice $\Theta$ is simply equal to $(2^{\aleph_0})^+$, it has very interesting behaviour in models of the Axiom of Determinacy.

Other examples of properties that can be used to obtain strong models of determinacy are ``all sets of reals are Suslin'' or ``all sets of reals are universally Baire.'' Being Suslin is a generalization of being analytic. More precisely, a set of reals is \emph{Suslin} if it is the projection of a tree on $\omega \times \kappa$ for some ordinal $\kappa$. Woodin and Steel determined the exact large cardinal strength of the theory ``$\AD$ + all sets of reals are Suslin'' \cite{St09DMT,St}:

\begin{theorem}[Steel, Woodin, \cite{St09DMT,St}] \label{thm:allsetsSuslin}
  The following theories are equiconsistent (over $\ZF$):
  \begin{enumerate}
      \item $\AD \, + \text{ all sets of reals are Suslin}$,
      \item $\ZFC \, + $ there is a cardinal $\lambda$ that is a limit of Woodin cardinals and a limit of ${<}\lambda$-strong cardinals.
  \end{enumerate}
  \end{theorem}

By results of Martin and Woodin, see \cite[Theorems 9.1 and 9.2]{St09DMT}, assuming $\AD$, the statement ``all sets of reals are Suslin'' is equivalent to the Axiom of Determinacy for games on reals ($\AD_\bR$).
Being universally Baire is a strengthening of being Suslin that was introduced implicitly in \cite{SchVa83} and explicitly by Feng, Magidor and Woodin \cite{FMW92}. 

\begin{definition}[Feng, Magidor, Woodin \cite{FMW92}] 
A subset $A$ of a topological space $Y$ is \emph{universally Baire} if $f^{-1}(A)$ has the property of Baire in any topological space $X$, where $f\colon X\rightarrow Y$ is continuous.
\end{definition} 

The exact consistency strength of the statement that all sets of reals are universally Baire under determinacy was conjectured by Sargsyan, in 2014, after he was able to obtain an upper bound with Larson and Wilson \cite{LSW} via an extension of Woodin's famous derived model theorem. One fact that makes their argument particularly interesting is that no model of the form $L(\Pot(\bR))$ is a model of $\text{``}\AD + $ all sets of reals are universally Baire.'' Universal Baireness is not only widely used across set theory but a crucial property in inner model theory: Universally Baire iteration strategies (canonically coded as a set of reals) can be extended from countable to uncountable iterations (see, for example, \cite{Sa13DIM}). The following theorem proves \textbf{Sargsyan's conjecture} by showing that the upper bound Larson, Sargsyan, and Wilson obtained is optimal:

\begin{theorem}[Larson, Sargsyan, Wilson, \cite{LSW}, Müller, \cite{Mu21}] \label{thm:allsetsuB}
The following theories are equiconsistent (over $\ZF$):
  \begin{enumerate}
      \item $\AD \, + \text{ all sets of reals are universally Baire}$,
      \item $\ZFC \, + $ there is a cardinal that is a limit of Woodin cardinals and a limit of strong cardinals.
  \end{enumerate}
  \end{theorem}

To construct and analyze the relevant models to prove the direction $Con(1.) \Rightarrow Con(2.)$ in this theorem, instead of just considering two hierarchies -- determinacy axioms and inner models with large cardinals -- a third hierarchy is used to reach higher levels in the other two. These three hierarchies together form what Steel calls the \textbf{triple helix} of inner model theory. The new hierarchy goes back to Woodin and Sargsyan and consists of canonical models called \textbf{hybrid mice}, or \textbf{hod mice}. These models are not only enhanced by large cardinals witnessed by extenders on their sequence, but also equipped with partial iteration strategies for themselves, see \cite{Sa15}. The strength of these models intuitively comes from the descriptive set theoretic complexity of these partial iteration strategies. 

The name \textbf{hod mouse} comes from the fact that these mice naturally occur as the result of analyses of $\HOD$, the hereditarily ordinal definable sets, in various models of determinacy. This analysis was pioneered by Steel and Woodin \cite{St95,StW16} 
in the model $L(\bR)$, as well as in $L[x][g]$ for a cone of reals $x$, where $g$ is generic for L\'{e}vy collapsing the least inaccessible cardinal in $L[x]$ (both under a determinacy hypothesis). It was extended to larger models of determinacy by Sargsyan, Trang, and others \cite{Sa15,Tr14,AtSa19,SaTrLSA}. In \cite{MuSa} we showed how to analyze $\HOD$ in $M_n(x)[g]$, for a cone of reals $x$, where $g$ is generic for L\'{e}vy collapsing the least inaccessible cardinal in $M_n(x)$ (under a determinacy hypothesis). 

The technical innovation behind the direction $Con(1.) \Rightarrow Con(2.)$ in Theorem \ref{thm:allsetsuB} is a new translation procedure to translate hybrid mice into mice with a limit of Woodin and strong cardinals \cite{Mu21}. This required an iterability proof for models obtained via a novel backgrounded construction. In \cite{Mu21} it is shown that the resulting models are countably iterably, meaning that countable substructures are iterable, and, in fact, a bit more. But the following natural question is left open:

\begin{question}
    Is there a translation procedure that yields fully iterable mice with a limit of Woodin and strong cardinals (when applied to suitable hybrid mice)?
\end{question}

%
%
%
\bibliographystyle{splncs04}
\bibliography{References}

\begin{thebibliography}{10}
\providecommand{\url}[1]{\texttt{#1}}
\providecommand{\urlprefix}{URL }
\providecommand{\doi}[1]{https://doi.org/#1}

\bibitem{AgMu}
Aguilera, J.P., Müller, S.: {Projective Games on the Reals}. Notre Dame
  Journal of Formal Logic  \textbf{61},  573--589 (2020).
  \doi{10.1215/00294527-2020-0027}

\bibitem{AgMu19}
Aguilera, J.P., Müller, S.: {The consistency strength of long projective
  determinacy}. {J. Symb. Log.}  \textbf{85}(1),  338–366 (2020).
  \doi{10.1017/jsl.2019.78}

\bibitem{AMS}
Aguilera, J.P., Müller, S., Schlicht, P.: Long games and {$\sigma$-projective}
  sets. Ann. Pure Appl. Log.  \textbf{172},  102939 (2021).
  \doi{10.1016/j.apal.2020.102939}

\bibitem{ANS01}
Andretta, A., Neeman, I., Steel, J.R.: The domestic levels of {$K^c$} are
  iterable. {Israel Journal of Mathematics}  \textbf{125},  157--201 (01 2001).
  \doi{10.1007/BF02773379}

\bibitem{AtSa19}
Atmai, R., Sargsyan, G.: {Hod up to $AD_{\mathbb{R}}$ + $\Theta$ is
  measurable}. Ann. Pure Appl. Log.  \textbf{170}(1),  95--108 (2019).
  \doi{10.1016/j.apal.2018.08.013}

\bibitem{CMM1}
Carroy, R., Medini, A., Müller, S.: {Every zero-dimensional homogeneous space
  is strongly homogeneous under determinacy}. J. Math. Log.  \textbf{20},
  2050015 (2020). \doi{10.1142/S0219061320500154}

\bibitem{CSS05}
Casacuberta, C., Scevenels, D., Smith, J.H.: {Implications of large-cardinal
  principles in homotopical localization}. Adv. Math.  \textbf{197}(1),
  120--139 (2005). \doi{10.1016/j.aim.2004.10.001}

\bibitem{DW87}
Dales, H.G., Woodin, W.H.: {An Introduction to Independence for Analysts}.
  Lond. Math. Soc. Lecture Note Ser., Cambridge University Press (1987).
  \doi{10.1017/CBO9780511662256}

\bibitem{Da64}
Davis, M.: {Infinite games of perfect information}. In: Dresher, M., Shapley,
  L.S., Tucker, A.W. (eds.) Advances in Game Theory, vol.~52, pp. 85--101
  (1964)

\bibitem{EM02}
Eklof, P.C., Mekler, A.H.: {Almost Free Modules}, North-Holland Mathematical
  Library, vol.~65. {North-Holland Publishing Co.} (2002)

\bibitem{Fa11}
Farah, I.: {All automorphisms of the Calkin algebra are inner}. Ann. Math.
  \textbf{173}(2),  619--661 (2011). \doi{10.4007/annals.2011.173.2.1}

\bibitem{FMW92}
Feng, Q., Magidor, M., Woodin, W.H.: {Universally Baire Sets of Reals}. In:
  Judah, H., Just, W., Woodin, W.H. (eds.) Set Theory of the Continuum.
  Mathematical Sciences Research Institute Publications, Springer (1992).
  \doi{10.1007/978-1-4613-9754-0\_15}

\bibitem{Fl82}
Fleissner, W.G.: {If all normal Moore spaces are metrizable, then there is an
  inner model with a measurable cardinal}. Trans. Amer. Math. Soc.
  \textbf{273},  365--373 (1982). \doi{10.1090/S0002-9947-1982-0664048-8}

\bibitem{GS53}
Gale, D., Stewart, F.M.: {Infinite games with perfect information}. In: Kuhn,
  H.W., Tucker, A.W. (eds.) Contrib. to the Theory of Games, vol.~2, pp.
  245--266 (1953)

\bibitem{Ha78}
Harrington, L.: {Analytic Determinacy and $0^\#$}. {J. Symb. Log.}
  \textbf{43},  685--693 (1978). \doi{10.2307/2273508}

\bibitem{Je72}
Jensen, R.B.: {The fine structure of the constructible hierarchy}. {Annals of
  Mathematical Logic}  \textbf{4},  229--308 (1972).
  \doi{10.1016/0003-4843(72)90001-0}

\bibitem{JSSS09}
Jensen, R.B., Schimmerling, E., Schindler, R., Steel, J.R.: Stacking mice. {J.
  Symb. Log.}  \textbf{74}(1),  315--335 (2009). \doi{10.2178/jsl/1231082314}

\bibitem{KS85}
Kechris, A., Solovay, R.: {On the Relative Consistency Strength of Determinacy
  Hypotheses}. {Trans. Amer. Math. Soc.}  \textbf{290}(1),  179--211 (1985).
  \doi{10.1090/S0002-9947-1985-0787961-2}

\bibitem{La}
Larson, P.B.: A brief history of determinacy. In: Kechris, A.S., Löwe, B.,
  Steel, J.R. (eds.) Large Cardinals, Determinacy and Other Topics, p. 3–60
  (2020). \doi{10.1017/9781316863534.002}

\bibitem{LaSa21}
Larson, P.B., Sargsyan, G.: {Failure of square in
  $\mathbb{P}_{\operatorname{max}}$ extensions of Chang models}  (2021)

\bibitem{LSW}
Larson, P.B., Sargsyan, G., Wilson, T.M.: {A model of the Axiom of Determinacy
  in which every set of reals is universally Baire}  (2018)

\bibitem{La76}
Laver, R.: {On the consistency of Borel's conjecture}. Acta Math.
  \textbf{137},  151--169 (1976). \doi{10.1007/BF02392416}

\bibitem{Ma70}
Martin, D.A.: {Measurable cardinals and analytic games}. {Fundamenta
  Mathematicae}  \textbf{66},  287--291 (1970). \doi{10.4064/fm-66-3-287-291}

\bibitem{Ma75}
Martin, D.A.: {Borel Determinacy}. {Ann. Math.}  \textbf{102}(2),  363--371
  (1975). \doi{10.2307/1971035}

\bibitem{MaSt89}
Martin, D.A., Steel, J.R.: {A Proof of Projective Determinacy}. {J. Amer. Math.
  Soc.}  \textbf{2}(1),  71--125 (1989). \doi{10.2307/1990913}

\bibitem{MS94}
Mitchell, W.J., Steel, J.R.: Fine structure and iteration trees, Lec. Not.
  Log., vol.~3. Springer-Verlag, Berlin, New York (1994)

\bibitem{Uh16}
Uhlenbrock~(now M\"uller), S.: {Pure and Hybrid Mice with Finitely Many Woodin
  Cardinals from Levels of Determinacy}. Ph.D. thesis, University of Münster
  (2016)

\bibitem{Mu19}
Müller, S.: {The axiom of determinacy implies dependent choice in mice}.
  {Mathematical Logic Quarterly}  \textbf{65}(3),  370--375 (2019).
  \doi{10.1002/malq.201800077}

\bibitem{Mu20rev}
Müller, S.: {Four papers on the large cardinal strength of
  $\operatorname{PFA}$ via core model induction}. Bull. Symb. Log.
  \textbf{26}(1),  89--92 (2020). \doi{10.1017/bsl.2020.6}

\bibitem{Mu21}
Müller, S.: The consistency strength of determinacy when all sets are
  universally {B}aire  (2021), submitted

\bibitem{MuSa}
Müller, S., Sargsyan, G.: H{OD} in inner models with {W}oodin cardinals. J.
  Symb. Log.  (2021). \doi{10.1017/jsl.2021.61}

\bibitem{MSW20}
Müller, S., Schindler, R., Woodin, W.: {Mice with Finitely many Woodin
  Cardinals from Optimal Determinacy Hypotheses}. J. Math. Log.  \textbf{20}
  (2020). \doi{10.1142/S0219061319500132}

\bibitem{Ne95}
Neeman, I.: {Optimal Proofs of Determinacy}. {Bull. Symb. Log.}  \textbf{1}(3),
   327--339 (09 1995). \doi{10.2307/421159}

\bibitem{Ne99}
Neeman, I.: {Games of Countable Length}. In: Cooper, S.B., Truss, J.K. (eds.)
  {Sets and Proofs}, p. 159–196. Lond. Math. Soc. Lecture Note Ser. (1999).
  \doi{10.1017/CBO9781107325944.009}

\bibitem{Ne02wlw}
Neeman, I.: {Inner models in the region of a Woodin limit of Woodin cardinals}.
  {Ann. Pure Appl. Log.}  \textbf{116}(1),  67 -- 155 (2002).
  \doi{10.1016/S0168-0072(01)00103-8}

\bibitem{Ne02}
Neeman, I.: {Optimal Proofs of Determinacy II}. {J. Math. Log.}  \textbf{2}(2),
   227--258 (11 2002). \doi{10.1142/S0219061302000175}

\bibitem{Ne04}
Neeman, I.: The Determinacy of Long Games, De Gruyter series in logic and its
  applications, vol.~7. De Gruyter (2004). \doi{10.1515/9783110200065}

\bibitem{Ny80}
Nyikos, P.J.: {A provisional solution to the normal Moore space problem}. Proc.
  Amer. Math. Soc.  \textbf{78},  429--435 (1980).
  \doi{10.1090/S0002-9939-1980-0553389-4}

\bibitem{Sa13DIM}
Sargsyan, G.: Descriptive inner model theory. The Bulletin of Symbolic Logic
  \textbf{19}(1),  1--55 (2013). \doi{10.2178/bsl.1901010}

\bibitem{Sa15}
Sargsyan, G.: {Hod Mice and the Mouse Set Conjecture}, Memoirs of the Amer.
  Math. Soc., vol.~236 (2015). \doi{10.1090/memo/1111}

\bibitem{Sa17}
Sargsyan, G.: {Translation procedures in descriptive inner model theory}. In:
  {Foundations of Mathematics}, vol.~690, pp. 205--223. Amer. Math. Soc.
  (2017). \doi{10.1090/conm/690/13869}

\bibitem{Sa21}
Sargsyan, G.: {Announcement of recent results in descriptive inner model
  theory}  (2021)

\bibitem{SaTrLSA}
Sargsyan, G., Trang, N.D.: {The Largest Suslin Axiom} (2016)

\bibitem{SaTr}
Sargsyan, G., Trang, N.D.: {The exact consistency strength of the generic
  absoluteness for the universally Baire sets}  (2019)

\bibitem{SaTrB}
Sargsyan, G., Trang, N.D.: {Sealing from iterability}. Trans. Amer. Math. Soc.
  Ser. B  \textbf{8},  229--248 (2021). \doi{10.1090/btran/65}

\bibitem{SaTrC}
Sargsyan, G., Trang, N.D.: {Sealing of the universally Baire sets}. Bull. Symb.
  Log.  (2021). \doi{10.1017/bsl.2021.29}

\bibitem{SchVa83}
Schilling, K., Vaught, R.: Borel games and the {B}aire property. Trans. Amer.
  Math. Soc.  \textbf{279}(1),  411--428 (1983). \doi{10.2307/1999393}

\bibitem{SchStZe02}
Schindler, R., Steel, J.R., Zeman, M.: {Deconstructing inner model theory}. {J.
  Symb. Log.}  \textbf{67},  721--736 (2002). \doi{10.2178/jsl/1190150106}

\bibitem{Sh74}
Shelah, S.: {Infinite abelian groups, Whitehead problem and some
  constructions}. Israel J. Math.  \textbf{18},  243--256 (1974).
  \doi{10.1007/BF02757281}

\bibitem{St88}
Steel, J.R.: {Long games}. In: Cabal Seminar 81-85, pp. 56--97. Lecture Notes
  in Math. 1333, Springer Verlag (1988). \doi{10.1007/BFb0084970}

\bibitem{St95}
Steel, J.R.: {$\operatorname{HOD}^{L(\mathbb{R})}$ is a Core Model below
  $\Theta$}. {Bull. Symb. Log.}  \textbf{1}(1),  75--84 (1995).
  \doi{10.2307/420947}

\bibitem{St}
Steel, J.R.: {An optimal consistency strength lower bound for
  $\mathsf{AD}_{\mathbb{R}}$}  (2008)

\bibitem{St08dermod}
Steel, J.R.: {Derived models associated to mice}. In: {Computational Prospects
  of Infinity - Part I: Tutorials}, vol.~14, pp. 105--193. World Scientific
  (2008). \doi{10.1142/9789812794055\_0003}

\bibitem{St09DMT}
Steel, J.R.: {The derived model theorem}. In: Cooper, S.B., Geuvers, H.,
  Pillay, A., Väänänen, J. (eds.) {Logic Colloquium 2006}, p. 280–327
  (2009). \doi{10.1017/CBO9780511605321.014}

\bibitem{St10}
Steel, J.R.: {An Outline of Inner Model Theory}, pp. 1595--1684. Springer
  (2010). \doi{10.1007/978-1-4020-5764-9\_20}

\bibitem{StW16}
Steel, J.R., Woodin, W.H.: {HOD as a core model}. In: Kechris, A.S., Löwe, B.,
  Steel, J.R. (eds.) {Ordinal definability and recursion theory}, pp. 257--346
  (2016). \doi{10.1017/CBO9781139519694.010}

\bibitem{Tr13}
Trang, N.D.: {Generalized Solovay Measures, the HOD Analysis, and the Core
  Model Induction}. Ph.D. thesis, University of California at Berkeley (2013)

\bibitem{Tr14}
Trang, N.D.: {$\HOD$ in natural models of $\AD^+$}. {Ann. Pure Appl. Log.}
  \textbf{165}(10),  1533 -- 1556 (2014). \doi{10.1016/j.apal.2014.04.006}

\bibitem{Wo55}
Wolfe, P.: The strict determinateness of certain infinite games. {Pacific
  Journal of Mathematics}  \textbf{5},  841--847 (1955).
  \doi{10.2140/pjm.1955.5.841}

\bibitem{Ze02}
Zeman, M.: {Inner Models and Large Cardinals}, De Gruyter series in logic and
  its applications, vol.~5. {De Gruyter} (2002). \doi{10.1515/9783110857818}

\bibitem{Zh15}
Zhu, Y.: {Realizing an $\AD^+$ model as a derived model of a premouse}. {Ann.
  Pure Appl. Log.}  \textbf{166},  1275--1364 (2015).
  \doi{10.1016/j.apal.2015.05.002}

\end{thebibliography}
%

\end{document}